\documentclass[11pt]{article}
\usepackage{amsmath,amsfonts,amssymb}
\usepackage{epsf}

\title{A strengthening and a multipartite generalization
       of the Alon-Boppana-Serre Theorem}

\author{Bojan Mohar\thanks{Supported in part by the
  Research Grant P1--0297 of ARRS (Slovenia), by an NSERC Discovery Grant (Canada)
  and by the Canada Research Chair program.}~\thanks{On leave from:
  IMFM \& FMF, Department of Mathematics, University of Ljubljana, Ljubljana,
  Slovenia.}\\
  {Department of Mathematics}\\
  {Simon Fraser University}\\
  {Burnaby, B.C. V5A 1S6}
  \\ email: {\tt mohar@sfu.ca}
}

\date{}

\newtheorem{theorem}{Theorem}[section]
\newtheorem{lemma}[theorem]{Lemma}
\newtheorem{corollary}[theorem]{Corollary}
\newtheorem{conjecture}[theorem]{Conjecture}

\newcommand{\DEF}[1]{{\em #1\/}}

\newcommand{\XX}{\ensuremath{\mathbb X}}
\newcommand{\TT}{\ensuremath{\mathbb T}}
\newcommand{\ZZ}{\ensuremath{\mathbb Z}}

\newcommand{\eopf}{\raisebox{0.8ex}{\framebox{}}}
\newcommand\slika[1]{\begin{center}\hskip 0.1mm\epsffile{#1}\end{center}}
\renewcommand\l{\lambda}
\newenvironment{proof}%
{\noindent{\bf Proof.}\ }%
{\hfill\eopf\par\bigskip}%

\begin{document}

\maketitle

\begin{abstract}
The Alon-Boppana theorem confirms that for every $\varepsilon>0$ and every
integer $d\ge3$, there are only finitely many $d$-regular graphs
whose second largest eigenvalue is at most $2\sqrt{d-1}-\varepsilon$.
Serre gave a strengthening showing that a positive proportion of eigenvalues
of any $d$-regular graph must be bigger than $2\sqrt{d-1}-\varepsilon$.
We provide a multipartite version of this result.
Our proofs are elementary and work also in the case when graphs are not regular.
In the simplest, monopartite case, our result extends the Alon-Boppana-Serre result
to non-regular graphs of minimum degree $d$ and bounded maximum degree.
The two-partite result shows that for every $\varepsilon>0$ and any positive
integers $d_1,d_2,d$, every $n$-vertex graph of maximum degree
at most $d$, whose vertex set is the union of (not necessarily disjoint) subsets
$V_1,V_2$, such that every vertex in $V_i$ has at least $d_i$
neighbors in $V_{3-i}$ for $i=1,2$, has $\Omega_\varepsilon(n)$ eigenvalues that are larger than $\sqrt{d_1-1}+\sqrt{d_2-1}-\varepsilon$.
Finally, we strengthen the Alon-Boppana-Serre theorem by showing that the lower
bound $2\sqrt{d-1}-\varepsilon$ can be replaced by $2\sqrt{d-1} + \delta$ for
some $\delta>0$ if graphs have bounded ``global girth''. On the other side of
the spectrum, if the odd girth is large, then we get an Alon-Boppana-Serre type
theorem for the negative eigenvalues as well.
\end{abstract}

\bigskip

{\bf Keywords:}
spectral radius, eigenvalue, Ramanujan graph, universal cover.

\medskip

{\bf Math.\,Subj.\,Classification:} {05C50}

\bigskip

\section{Introduction}

After the breakthrough paper by Alon and Milman \cite{AM} in 1985,
it became apparent that regular graphs, whose \DEF{spectral gap}
(i.e.~the difference between the largest and the second largest eigenvalue)
is large, posses some extraordinary properties, like unusually fast expansion
and resemblance to random graphs. This led to the definition of
\DEF{Ramanujan graphs}. These are $d$-regular graphs whose second largest
eigenvalue does not exceed the value $2\sqrt{d-1}$.
Lubotzky, Phillips, and Sarnak \cite{LPS},
and independently Margulis \cite{Ma}, were the first to show that Ramanujan graphs
exist. Their constructions are based on number theory and work when
the degree $d$ is equal to $p+1$ for some prime $p$. Later, several new
constructions were discovered, showing that Ramanujan graphs exist for every
degree $d\ge3$ which is of the form $p^k+1$ for some integer $k\ge1$ and some
prime $p$, see Morgenstern \cite{Morg}.

It is not immediately clear why the special choice of $2\sqrt{d-1}$ is taken
when defining Ramanujan graphs. One reason is that this is the spectral radius
of the infinite $d$-regular tree, which is the universal cover for all
$d$-regular graphs. Another reason is the following result
of Alon and Boppana (see~\cite{Al}) which shows that this is the smallest
number that makes sense.

\begin{theorem}[Alon-Boppana]
\label{thm:AB}
For every $d\ge2$ and every $\varepsilon>0$, there are only finitely many
$d$-regular graphs whose second largest eigenvalue is at most
$2\sqrt{d-1}-\varepsilon$.
\end{theorem}

Alternative proofs of Theorem \ref{thm:AB} were given by Friedman \cite{Fr}
and by Nilli \cite{Ni1}, who has recently further simplified her arguments
in \cite{Ni2}. Actually, the proofs in \cite{Fr,Ni2} imply a stronger version
of Theorem \ref{thm:AB} by making the same conclusion for more
eigenvalues than just the second largest one.
This strengthening is attributed to Serre \cite{Se} (see also \cite{DSV,FL,HLW}),
who stated the following quantitative version of the Alon-Boppana Theorem:

\begin{theorem}[Serre]
\label{thm:Serre}
For every positive integer $d$ and every $\varepsilon>0$, there exists a constant
$c=c(d,\varepsilon)$ such that every $d$-regular graph of order $n$
has at least $cn$ eigenvalues that are larger than $2\sqrt{d-1}-\varepsilon$.
\end{theorem}




In this paper we give a multipartite generalization of the Alon-Boppana
Theorem, see Theorem \ref{thm:main}. The Ramanujan value $2\sqrt{d-1}$
is replaced by the spectral radius of the universal covering
tree of the multipartite parameters (cf.\ Section~\ref{sect:2} for definitions).
Our proof has similarities with
Nilli's proof \cite{Ni2}, and seems to be even simpler if restricted to the
special case of $d$-regular graphs. The main step is
based on the interlacing theorem and is entirely elementary.

Our proofs work also in the case when graphs are not regular.
In the simplest, monopartite case, our result extends the Alon-Boppana-Serre result
to non-regular graphs of minimum degree $d$ and bounded maximum degree.
A strengthening of this form has been obtained previously by Hoory \cite{Ho}.
In the next simplest two-partite case it is shown that for every $\varepsilon>0$ and any positive
integers $d_1,d_2,d$, every $n$-vertex graph of maximum degree
at most $d$, whose vertex set is the union of (not necessarily disjoint) subsets
$V_1,V_2$, such that every vertex in $V_i$ has at least $d_i$
neighbors in $V_{3-i}$ for $i=1,2$, has $\Omega_\varepsilon(n)$ eigenvalues that are larger than $\sqrt{d_1-1}+\sqrt{d_2-1}-\varepsilon$.

After submission of this paper, S. Cioab\v a informed us about some related
work. Greenberg \cite{Gr} obtained a generalized version of the Serre theorem in a similar form as ours, but only claiming that there are eigenvalues whose absolute value is larger than $2\sqrt{d-1}-\varepsilon$. Cioab\v a \cite{Ci}
improved Greenberg's work to the same form as given in Theorem \ref{thm:Serre}. Greenberg's result also appears in \cite{LN}.

In the last section we tailor the proofs to obtain a strengthening
of the Alon-Boppana-Serre theorem by showing that the lower
bound $2\sqrt{d-1}-\varepsilon$ can be replaced by $2\sqrt{d-1} + \delta$ for
some $\delta>0$ if graphs have bounded universal girth (see Section \ref{sect:6}
for the definition). On the other side of the spectrum, if the \DEF{odd girth}
(i.e.\ the length of a shortest odd cycle) is large, then we obtain
an Alon-Boppana-Serre type theorem for the negative eigenvalues.

If $G$ is a (finite) graph, we denote by $\l_i=\l_i(G)$ the $i$th largest eigenvalue
of the adjacency matrix $A(G)$ of $G$, respecting multiplicities.
The largest eigenvalue of $G$, $\rho(G)=\l_1(G)$, is also referred to
as the \DEF{spectral radius} of $G$. It follows from the Perron-Frobenius
theorem (see, e.g.~\cite{HJ}) that $\rho(G)$ is an eigenvalue of $G$ that
has an eigenvector $x$ whose coordinates are all non-negative. Moreover,
if $G$ is connected, then $x$ is strictly positive.

If $r\ge 1$ is an integer, a set $S$ of vertices of a graph $G$ is said to be
\DEF{$r$-apart} if any two vertices in $S$ are at distance at least $r+1$ in $G$.
We denote by $\alpha_r(G)$ the maximum cardinality of a vertex set in $G$ that
is $r$-apart. Note that $\alpha_1(G)$ is the usual independence number of the
graph.

Let $G$ be a graph, $v\in V(G)$, and let $r$ be an integer.
We denote by $G_r(v)$ the induced subgraph of $G$ on vertices that are
at distance at most $r$ from $v$. The subgraph $G_r(v)$ is called the
\DEF{$r$-ball} around $v$ in $G$.

We allow infinite graphs, but they will always be locally finite.
In particular, the $r$-ball around any vertex of a graph $G$ is always
finite.

\section{Universal covers and subcovers}

\label{sect:2}

Let $D=[d_{ij}]_{i,j=1}^t$ be a square matrix of order $t\ge1$, whose entries
$d_{ij}$ are non-negative integers. For $i=1,\dots,t$, we define the
\DEF{$i^{th}$ degree} in $D$ as the integer $d_i = \sum_{j=1}^t d_{ij}$.
Suppose that $D$ further satisfies the following conditions:
\begin{itemize}
\item[\bf (D1)] If $d_{ij}=0$, then also $d_{ji}=0$.
\item[\bf (D2)] The graph of $D$ is connected, i.e., for every
$i,k\in \{1,\dots,t\}$ there are integers $m_1,m_2,\dots,m_s$ in $\{1,\dots,t\}$,
where $m_1=i$, $m_s=k$, and $d_{m_j m_{j+1}} > 0$ for $j=1,\dots,s-1$.
\item[\bf (D3)] For every sequence of (distinct) integers $m_1,m_2,\dots,m_s$ in
$\{1,\dots,t\}$, we have
$$
    d_{m_1 m_2}d_{m_2 m_3}\cdots d_{m_{s-1}m_s}d_{m_s m_1} =
    d_{m_1 m_s}d_{m_s m_{s-1}}\cdots d_{m_3 m_2}d_{m_2 m_1}.
$$
\end{itemize}
Such a matrix is called a \DEF{$t$-partite degree matrix}.

Let $D=[d_{ij}]$ be a $t$-partite degree matrix. If a graph $G$ admits a
partition of its vertex set into $t$ classes, $V(G)=U_1\cup\cdots\cup U_t$,
such that every vertex in $U_i$ has precisely $d_{ij}$ neighbors in $U_j$,
for all $i,j=1,\dots,t$, then we say that $D$ is a \DEF{$t$-partite degree matrix}
for $G$. The corresponding partition $U_1\cup\cdots\cup U_t$ is said to be an
\DEF{equitable partition} for $D$; see, e.g.~\cite{GR}.

\begin{lemma}
\label{lem:universal cover}
Let\/ $D$ be a $t$-partite degree matrix.

{\rm (a)} There exists a finite graph $G$ whose degree matrix is $D$.

{\rm (b)} There exists a tree $T_D$ whose degree matrix is $D$.
The tree $T_D$ is determined up to isomorphism.
\end{lemma}

\begin{proof}
(a) First we remark that the condition (D3) implies that the set of equalities
$n_i d_{ij} = n_j d_{ji}$, $i,j\in \{1,\dots,t\}$, has a positive solution
$n_1,\dots,n_t$. Since all $d_{ij}$ are integers, there is a solution whose
values $n_i$ ($i=1,\dots,t$) are positive integers. To obtain a graph $G$,
we take vertex sets $U_i$ of cardinalities $n_i$ for $i=1,\dots,t$, and
join $U_i$ and $U_j$ so that the edges between them form
a $(d_{ij},d_{ji})$-biregular bipartite graph. Then it is clear that
$D$ is a $t$-partite degree matrix for $G$.

(b) To get $T_D$, we just take what is known as the universal cover of
the graph $G$ obtained in part (a).

We add a short proof of existence of $T_D$ that does not use
the property (D3) which is needed in (a).
Let us first assume that $d_i\ge2$ for $i=1,\dots,t$. This case is intuitively
clear and we leave the details of the proof for the reader. Note that $T_D$
is always infinite in this case.

The rest of the proof is by induction on $t$. We may assume that $d_t\le1$.
If $t=1$, then $T_D$ is either a single vertex (if $d_1=0$) or an edge
(if $d_1=1$). If $t>1$, then (D2) implies that $d_t=1$ and the non-zero element
in row $t$ of $D$ is not $d_{tt}$. Thus, there is a unique $j<t$ such that
$d_{tj}=1$. Let $p=d_{jt}$. By (D1), we conclude that $p>0$ and all other
elements in the column $t$ of $D$ are zero. Let $D'$ be the submatrix of $D$
obtained by deleting the last row and the last column. Since this operation
acts like removing a vertex of degree 1 from a graph, $D'$ still satisfies
(D1)--(D2), and hence we can apply the induction hypothesis to find the tree
$T_{D'}$. Finally, we obtain $T_D$ by adding, to each vertex in $V_j$,
$p$ pendant edges. All new vertices are of degree 1 and form the class $V_t$
in $T_D$.
\end{proof}

The tree $T_D$ is called the \DEF{universal cover} of the multipartite
degree matrix $D$. Let $V_1\cup\cdots\cup V_t$ be the corresponding equitable
partition of $V(T_D)$. If $G$ is any graph whose $t$-partite degree matrix is
$D$, there is a covering projection $\pi_G : T_D\to G$ which maps vertices in
$V_i$ onto the $i$th class of the equitable partition of $V(G)$.

Covering projections, universal covers and equitable partitions are regularly
used in algebraic graph theory. In the sequel we shall introduce a weaker notion,
distinguished by the prefix `sub', in which only those properties that are
important for our main results will be preserved.

A degree matrix $D$ is said to be a \DEF{$t$-partite subdegree matrix} for
a graph $G$ if there is a graph homomorphism $\pi^D_G: T_D\to G$ which is
\DEF{locally $1$-$1$}, i.e., for each vertex $v\in V(T_D)$,
$\pi^D_G$ maps edges incident with $v$ injectively to the edges incident
with $\pi^D_G(v)$. The homomorphism $\pi^D_G: T_D\to G$ is called
a \DEF{subuniversal projection} and the tree $T_D$ is a \DEF{subuniversal cover}
of~$G$.

If $\pi^D_G: T_D\to G$ is a subuniversal projection, let
$U_i = \pi^D_G(V_i)\subseteq V(G)$, $i=1,\dots,t$. Then it is easy to see
that the (not necessarily disjoint) vertex-sets $U_1,\dots,U_t$ satisfy
the following condition: Every vertex in $U_i$ has at least $d_{ij}$ neighbors
in $U_j$, for all $i,j=1,\dots,t$.
This gives a necessary condition for existence of a subuniversal projection.
Unfortunately, this condition is not sufficient. But if we ask that
every vertex in $U_i$ has at least $d_{ij}+1$ neighbors in $U_j$, for all
$i,j=1,\dots,t$, then the existence of a subuniversal projection to
$G$ is easily verified.

\begin{theorem}
\label{thm:ball and tree}
Suppose that\/ $D$ is a subdegree matrix for a (possibly infinite) graph $G$,
and let\/ $T_D$ be the corresponding subuniversal cover.
Let $v\in V(G)$ and let $s\in V(T_D)$ be a vertex that is mapped to $v$ via
a subuniversal projection $\pi^D_G$. Then for every $r\ge0$, the spectral radius
of the $r$-ball in $G$ is at least as large as the spectral radius of
the corresponding $r$-ball in $T_D$, $\rho(G_r(v))\ge\rho(T_{D,r}(s))$.
\end{theorem}

\begin{proof}
The spectral radius of a connected graph $H$ can be expressed as
\begin{equation}
    \rho(H) = \limsup_{q\to \infty} (w_{2q}(H,u))^{1/(2q)},
    \label{eq:walks1}
\end{equation}
where $w_{2q}(H,u)$ denotes the number of closed walks of length $2q$ in $H$
starting at the vertex $u\in V(H)$.
Every closed walk in $T_{D,r}(s)$
starting at $s$ is projected by $\pi^D_G$ to a closed walk in $G_r(v)$
starting at $v$. The projection of these walks is 1-1, since $\pi^D_G$ is
locally 1-1. Hence,
\begin{equation}
    w_{2q}(G_r(v),v) \ge w_{2q}(T_{D,r}(s),s).
    \label{eq:walks2}
\end{equation}
This inequality in combination with (\ref{eq:walks1}) implies that
$\rho(G_r(v))\ge\rho(T_{D,r}(s))$.
\end{proof}

\section{The spectral radius of infinite trees}

If $G$ is a connected infinite (locally finite) graph, we define its
spectral radius $\rho(G)$ as
\begin{equation}
 \rho(G) = lim_{r\to\infty} \rho(G_r(v))
\label{eq:rho}
\end{equation}
where $v$ is any vertex of $G$. It is easy to see that the limit
exists (it may be infinite if the degrees of $G$ have no finite upper bound)
and that it is independent of the choice of $v$. The spectral radius of
infinite graphs defined above coincides with the notion obtained through the
spectral theory of linear operators in Hilbert spaces; we refer to
\cite{MW} for an overview.

The monotonicity property of the spectral radius of finite graphs implies that
for every connected finite graph $H$ and any proper subgraph $H'$ of $H$,
we have $\rho(H') < \rho(H)$.
Since $G$ is connected, infinite, and locally finite,
$G_r(v) \ne G_{r+1}(v) \ne G$ for every $r\ge0$.
This implies that
$$
   \rho(G_r(v)) < \rho(G_{r+1}(v)) < \rho(G).
$$

Let us remark that the spectral radius of an infinite $d$-regular
tree is equal to $2\sqrt{d-1}$, the value that appears in the definition of
Ramanujan graphs. This was proved by Kesten \cite{Ke}, see also
Dynkin and Malyutov \cite{DM}, Cartier \cite{Ca}, and Woess \cite{Wo}.
We will use the spectral radius of universal cover trees introduced
in the previous section to replace the Ramanujan bound $2\sqrt{d-1}$
with the corresponding bound suitable for our multipartite generalization.


In the special case when the graph is the infinite $d$-regular tree,
which shall be denoted by $\TT_d$, it is easy to determine the precise rate
of convergence in (\ref{eq:rho}).

\begin{theorem}
\label{thm:convergence}
For every integer $d\ge 2$, we have
$$\rho(\TT_{d,r}) > 2\sqrt{d-1}\,\Bigl(1 - \frac{\pi^2}{r^2} + O(r^{-3})\Bigr).$$
\end{theorem}

\begin{proof}
Let $w_q(G,v_0)$ denote the number of closed walks of length $q$. It will be convenient to consider the subtree $\TT'_d$ of $\TT_d$ which is equal to the connected component containing the vertex $v_0$ of the subgraph obtained after deleting an edge of $\TT_d$ incident with $v_0$. The vertex $v_0$ has degree $d-1$ in $\TT'_d$, while all other vertices still have degree $d$.
The tree $\TT'_d$ has a natural projection onto the one-way-infinite path $P_\infty$ (whose vertices we denote by the non-negative integers $0, 1, 2, \dots$) such that all vertices at distance $i$ from $v_0$ are mapped onto the vertex $i$ in $P_\infty$. Every closed walk (based at $v_0$) of length $2q$ in $\TT'_d$ is projected onto a closed walk in $P_\infty$ based at the vertex $0$. Moreover, the $r$-ball $\TT'_{d,r}(v_0)$
in $\TT'_d$ is projected onto the path $P_{r+1} \subset P_\infty$ on vertices
$0,1,\dots,r$.

Whenever we walk away from $v_0$ in $\TT'_d$, we have $d-1$ choices to do so. This implies that
\begin{equation}
    w_{2q}(\TT'_{d,r},v_0) = (d-1)^q \, w_{2q}(P_{r+1},0).
    \label{eq:walks T_d,r}
\end{equation}
When $q\to \infty$, the quantities raised to the power $1/(2q)$ tend to the spectral radii of the corresponding graphs, and we conclude that
$\rho(\TT'_{d,r})=\sqrt{d-1}\,\rho(P_{r+1}) = 2\sqrt{d-1}\,\cos(\frac{\pi}{r+2}) = 2\sqrt{d-1}\,\bigl(1 - \tfrac{\pi^2}{r^2} + O(r^{-3})\bigr)$.
Since $\TT'_{d,r}$ is a proper finite subgraph of $\TT_{d,r}$, this implies the (strict) inequality of the theorem.
\end{proof}

The rate of convergence is likely the same for more general universal covers of finite graphs. We propose the following conjecture.

\begin{conjecture}
\label{conj:convergence}
For every multipartite degree matrix $D$, there exists a constant $c=c(D)$
such that for every $s\in V(T_D)$, we have
$$\rho(T_{D,r}(s)) \ge \rho(T_D) - cr^{-2}.$$
\end{conjecture}

\section{Multipartite Ramanujan graphs}

In this section we introduce a generalized notion of Ramanujan graphs.
The following lemma shows that we cannot simply compare $\l_2(G)$ with
$\rho(T_D)$ as is the case for $d$-regular graphs.

\begin{lemma}
\label{lem:contains eigenvalues}
If\/ $D$ is a $t$-partite degree matrix, then all eigenvalues of\/ $D$ are
real and their algebraic multiplicity is equal to their geometric multiplicity.
If\/ $D$ is a multipartite degree matrix of a finite graph\/ $G$, then every
eigenvalue of $D$ is also an eigenvalue of\/ $G$. Moreover, $\rho(G)=\rho(D)$.
\end{lemma}

\begin{proof}
Let $n_1,\dots,n_t$ be a positive solution of the system
$n_i d_{ij} = n_j d_{ji}$, $i,j\in \{1,\dots,t\}$, which was shown to exist in
the proof of Lemma \ref{lem:universal cover}(a). If $R$ is the diagonal matrix
of order $t$ whose entry $R_{ii}$ is equal to $n_i^{1/2}$ ($i=1,\dots,t$),
then $RDR^{-1}$ is a symmetric matrix. This implies the first part of the lemma.

To verify the second part, let $\l$ be an eigenvalue of $D$, and
let $y=(y_i\mid i=1,\dots,t)$ be an eigenvector
for $\l$. Let $V_1\cup\cdots\cup V_t$ be the partition of $V(G)$ corresponding
to the degree matrix $D$. If we set $x_v = y_i$ for every $v\in V_i$, then it is
easy to see that $x=(x_v\mid v\in V(G))$ is an eigenvector of the adjacency
matrix of $G$ for the eigenvalue $\l$.

To prove the last claim, observe that the eigenvalue $\rho(D)$ has a positive eigenvector by the Perron-Frobenius Theorem. Its lift in $G$ is a positive
eigenvector of $G$ for the eigenvalue $\rho(D)$. Again, by applying
the Perron-Frobenius Theorem, we conclude that this eigenvector corresponds to the
largest eigenvalue of $G$.
\end{proof}

Let $D$ be a $t$-partite degree matrix. Let $k$ be the largest integer such
that $\l_k(D) \ge \rho(T_D)$. Note that $k$ exists since $\l_1(D) \ge \rho(T_D)$.
We say that a finite graph $G$ with degree matrix (resp.~subdegree matrix) $D$
is \DEF{$D$-Ramanujan} (resp.~\DEF{$D^+$-Ramanujan}) if
$\lambda_{k+1}(G)\le \rho(T_D)$. We believe that there is an abundance of
generalized Ramanujan graphs and propose the following conjectures
(in which we assume that the minimum degree of $D$ is at least 2).

\begin{conjecture}
\label{conj:exist D-Ramanujan2}
If there exists a $D$-Ramanujan graph for a multipartite degree matrix $D$,
then there exist infinitely many $D$-Ramanujan graphs.
\end{conjecture}

\begin{conjecture}
\label{conj:exist D-Ramanujan+}
If there exists a $D^+$-Ramanujan graph for a multipartite degree matrix $D$,
then there exist infinitely many $D^+$-Ramanujan graphs.
\end{conjecture}

\begin{conjecture}
\label{conj:exist D-Ramanujan}
If\/ $D$ is a degree matrix of order $t\ge2$, and $\l_2(D) < \rho(T_D)$,
then there exist infinitely many $D$-Ramanujan graphs.
\end{conjecture}

One cannot exclude the possibility that there exist infinitely many $D$-Ramanujan graphs for every degree matrix $D$, but our knowledge is too limited at this point to propose this as a conjecture.

\section{A generalized Alon-Boppana-Serre theorem}

\begin{theorem}
\label{thm:main}
Let\/ $D$ be a multipartite degree matrix, and let\/ $\rho_D = \rho(T_D)$.

{\rm (a)}
For every $\varepsilon>0$, there exists an integer $r=r(D,\varepsilon)$ such that
for every integer $k\ge1$ and for every graph $G$, if $D$ is a subdegree matrix
of $G$ and $\alpha_{2r+1}(G)\ge k$, then
$\l_k(G) \ge \rho_D - \varepsilon$.

{\rm (b)}
For $\varepsilon>0$ and every positive integer $\Delta$, there exists a constant
$c=c(D,\Delta,\varepsilon)>0$ such that every graph $G$ of order $n$, of maximum
degree at most $\Delta$ and with subdegree matrix $D$, has
at least $cn$ eigenvalues that are larger than $\rho_D-\varepsilon$.
\end{theorem}

\begin{proof}
(a) Let $r=r(D,\varepsilon)$ be the smallest integer such that
$\rho(T_{D,r})\ge \rho_D - \varepsilon$, and let $k$ and $G$ be as specified.
Since $\alpha_{2r+1}(G) \ge k$, there are vertices $v_1,\dots,v_k$ that are
$(2r+2)$-apart. The $r$-balls $G_r(v_1),\dots,G_r(v_k)$ around
these vertices are not only pairwise disjoint, but also form an induced subgraph
of $G$. By the eigenvalue interlacing property for induced subgraphs, we know that
$$
    \l_k(G) \ge \l_k(G_r(v_1)\cup\cdots\cup G_r(v_k)) \ge
    \min \{ \rho(G_r(v_i))\mid 1\le i\le k \}.
$$
By Theorem \ref{thm:ball and tree} and by our choice of $r$, we have
$$
   \rho(G_r(v_i)) \ge \rho(T_{D,r}(s_i)) \ge \rho_D - \varepsilon.
$$
This completes the proof of (a).

(b) This part follows from (a). It is just to be noted that
any $(2r+1)$-ball in $G$ contains at most
$B = \tfrac{\Delta}{\Delta-2}(\Delta-1)^{2r+1}$ vertices.
Thus, $\alpha_{2r+1}(G)\ge n/B$, and hence part (a) applies
with $c=B^{-1}$.
\end{proof}

It is worth mentioning that the condition involving $\alpha_{2r+1}(G)$ in
Theorem \ref{thm:main}(a) is necessary if we only assume that $D$ is a
subdegree matrix. Simple examples showing this are provided by the
family of all complete graphs $K_n$ whose second largest eigenvalue
is always equal to $-1$, or by the family of all complete bipartite
graphs $K_{m,n}$ whose second eigenvalue is 0.

For the special case when $D=[d]$, Theorem \ref{thm:convergence} gives the precise
description for the values $r(d,\varepsilon)$ and $c(d,\Delta,\varepsilon)$
in Theorem \ref{thm:main}. By Theorem~\ref{thm:convergence},
\begin{equation}
  r = r(d,\varepsilon) = \pi\biggl(\frac{2\sqrt{d-1}}{\varepsilon}\, \biggr)^{1/2}\,
                     \bigl(1 + O\bigl(d^{-1/4}\varepsilon^{1/2}\bigr) \bigr)
  \nonumber
\end{equation}
and
\begin{equation}
  c(d,\Delta,\varepsilon) = \tfrac{\Delta}{\Delta-2}(\Delta-1)^{-(2r+1)}
  \nonumber
\end{equation}
will do the job.

As an example, let us consider the following special case.
The bipartite degree matrix
\begin{equation}
   D = \left[\begin{matrix}
            0 & d_1 \\
            d_2 & 0
         \end{matrix}\right]
\label{eq:D}
\end{equation}
involves, in particular, all bipartite graphs with bipartition
$V= A\cup B$, whose degrees in $A$ are at least $d_1$ and whose degrees in
$B$ are at least $d_2$. The spectral radius of $T_D$ is (cf.~\cite{MW})
$$
   \rho(T_D) = \sqrt{d_1-1} + \sqrt{d_2-1}\,.
$$
Thus, only finitely many bipartite $(d_1,d_2)$-biregular graphs have their
$k$th eigenvalue ($k\ge2$) smaller than $\rho(T_D) - \varepsilon$.
Theorem \ref{thm:main} suggests the following strengthening, which we will prove
directly by using Theorem \ref{thm:ball and tree}.

\begin{corollary}
\label{cor:2partite}
Let $d_1\le d_2\le d$ be positive integers, and let ${\cal G}_{d_1,d_2}^d$ be
the set of all graphs whose maximum vertex degree is at most $d$ and whose
vertex set is the union of (not necessarily disjoint) subsets
$U_1,U_2$, such that every vertex in $U_i$ has at least $d_i$
neighbors in $U_{3-i}$ for $i=1,2$.
For every $\varepsilon>0$, every $n$-vertex graph $G\in {\cal G}_{d_1,d_2}^d$
has\/ $\Omega_\varepsilon(n)$ eigenvalues larger than
$\sqrt{d_1-1} + \sqrt{d_2-1} - \varepsilon$.
\end{corollary}

\begin{proof}
%
We claim that there exists a subuniversal projection $\pi^D_{G}$, where $D$ is
the degree matrix given in (\ref{eq:D}). The tree $T_D$ is $(d_1,d_2)$-biregular.
We map a vertex $v_0$ of degree $d_1$ in $T_D$ onto the vertex $u\in U_1$.
After fixing $v_0$, we extend the mapping to a locally 1-1
homomorphism in a greedy fashion (by taking the breadth-first search order
of vertices of $T_D$ starting at $v_0$) so that vertices of degree $d_i$ are mapped
to $U_i$, $i=1,2$.

Let $v_1$ be a neighbor of $v_0$ in $T_D$ and let $v = \pi_G^D(v_1)\in U_2$.
Theorem \ref{thm:ball and tree} shows that
for large enough $r=r(d_1,d_2,\varepsilon)$,
\begin{eqnarray}
  \rho(G_r(u)) &\ge& \rho(T_{D,r}(v_0))\ \ge\ \rho(T_D)-\varepsilon,
\label{eq:reachable}\\
  \rho(G_r(v)) &\ge& \rho(T_{D,r}(v_1))\ \ge\ \rho(T_D)-\varepsilon.
\label{eq:reachable2}
\end{eqnarray}
Since the maximum degree of $G$ is bounded by $d$, the $(2r+1)$-balls in $G$ have
bounded number of vertices, say at most $B$. Therefore,
$\alpha_{2r+1}(G)\ge n/B$, and so there are at least this many pairwise
non-adjacent induced $r$-balls around vertices in $G$. As before, the eigenvalue
interlacing theorem and (\ref{eq:reachable})--(\ref{eq:reachable2})
imply that linearly many eigenvalues
of $G$ are larger than

\medskip

\hskip 2.3cm
$\rho(T_D)-\varepsilon = \sqrt{d_1-1} + \sqrt{d_2-1} - \varepsilon$.
\end{proof}

\section{Global girth and Ramanujan graphs}

\label{sect:6}

All known Ramanujan graphs are Cayley graphs and their girth increases with their
order. We shall use the method of this paper to explain why the girth cannot be
bounded. Actually, we shall prove that a small girth condition implies that
$d$-regular graphs are ``far from being Ramanujan;''
see Theorem \ref{thm:girth} below.

Let $G$ be a graph. A closed walk $v_1v_2\dots v_kv_1$ is
\DEF{retracting-free} if $v_{i-1}\ne v_{i+1}$ for $i=1,\dots,k$
(where $v_0=v_k$ and $v_{k+1}=v_1$).
It is easy to see that if $G$ is a finite graph with minimum degree at least 2, then
for every vertex $v$ of $G$ there exists a retracting-free closed walk through $v$.

Let $g(v)$ be the length of a shortest retracting-free closed walk through $v$.
The \DEF{universal girth} of $G$, denoted by $m(G)$, is the smallest integer
$k$ such that every vertex in $G$ has a
retracting-free closed walk of length $k$. Let us observe that $m(G)$ is at most
the least common multiple of the values $g(v)$, $v\in V(G)$.
Also, if $G$ is vertex-transitive, then $m(G)$ is equal to the girth of~$G$.

\begin{figure}[htb]
\epsfxsize=6.4truecm
\slika{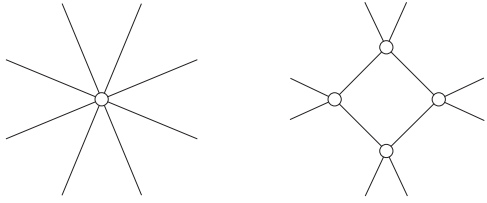}
\caption{Expanding a vertex into a cycle}
\label{fig:1}
\end{figure}

Let $\XX_{d,g}$ be the graph obtained from the $(d-2)g$-regular tree $T$ by
expanding each vertex $v\in V(T)$ into the cycle $C_v$ of length $g$, such that each
vertex of $C_v$ is incident with $d-2$ of the edges of $T$ incident with $v$.
See Figure \ref{fig:1} showing the case of $d=4$ and $g=4$.
The graph $\XX_{d,g}$ is the Cayley graph of the free product of $d-2$ copies
of $\ZZ_2$ and one copy of $\ZZ_g$ (with the natural generating set).

Paschke \cite{Pa} determined the spectral radius of $\XX_{d,g}$:

\begin{theorem}[Paschke]
\label{thm:Paschke}
For $d\ge3$ and\/ $g\ge3$, the graph\/ $\XX_{d,g}$ has spectral radius
$$
   \min_{s>0} \,(d-2)\, \phi\biggl(\frac{1+\cosh sg}{\sinh sg\sinh s}\biggr)
              + 2\cosh s > 2\sqrt{d-1},
$$
where $\phi(t) = \frac{\sqrt{1+t^2} - 1}{t}$.
\end{theorem}

Paschke~\cite{Pa} used this result to provide a non-trivial lower bound on
the spectral radius of infinite vertex-transitive graphs of the
given girth~$g$. He showed that a vertex transitive $d$-regular graph
containing a $g$-cycle has spectral radius at least $\rho(\XX_{d,g}$).
The formula in Theorem \ref{thm:Paschke} gives a lower bound of the form
$$
   2\sqrt{d-1} + \frac{2(d-2)}{(d-1)^{(g+1)/2}}\, h(d,g),
$$
where $h$ is a function such that such that for every $g\ge 3$,
$\lim_{d\to\infty} h(d,g)=1$, and for every $d\ge3$,
$\lim_{g\to\infty} h(d,g)=1$.

Now, we strengthen the Alon-Boppana-Serre theorem by showing that the lower bound
$2\sqrt{d-1}-\varepsilon$ can be replaced by
$2\sqrt{d-1} + \delta$ for some $\delta>0$ if graphs have bounded universal girth.

\begin{theorem}
\label{thm:girth}
For every $\Delta\ge d\ge3$ and every $g\ge3$, there exist $\delta>0$ and $c>0$
such that every $n$-vertex graph $G$ with minimum degree at least $d$, maximum
degree at most $\Delta$ and
universal girth at most $g$ has at least $\lceil cn\rceil$ eigenvalues that
are larger than $2\sqrt{d-1} + \delta$.
\end{theorem}

\begin{proof}
The proof follows the same pattern as the proof of Theorem \ref{thm:main}, except that we use the graph $\XX_{d,m}$, where $m=m(G)\le g$ is the universal girth of $G$,
playing the role of the universal cover $\TT_d$.
Here, we have to take $r$ large enough so that
$\rho((\XX_{d,m})_r) \ge 2\sqrt{d-1} + \delta$.
Such an $r$ exists because of (\ref{eq:rho}) and since
$\rho(\XX_{d,m}) > \rho(\TT_d) = 2\sqrt{d-1}$.
\end{proof}

It is straightforward to generalize the proof of Theorem \ref{thm:girth} to
the setting of degree matrices. What we need is just an analogue of the Paschke
theorem. However, we do not intend to dig into the details in this note.

\section{The other side of the spectrum}

As shown in the previous section, small universal girth yields improved
lower bounds on large eigenvalues, so Ramanujan graphs must have growing
girth. On the other hand, large girth has some further cosequences.
In particular, it shows that the negative eigenvalues satisfy the
Alon-Boppana-Serre property as well.

Let us first formulate the monopartite version for the negative eigenvalues.
It involves the notion of the \DEF{odd girth} of the graph, meaning the length
of a shortest cycle of odd length in the graph. (If $G$ is bipartite, then the
odd girth is $\infty$.) This result was obtained earlier by Friedman \cite{Fr}
and Nilli \cite{Ni2}; it also appears in Ciaob\v a \cite{Ci} (with a slightly weaker estimate of $\varepsilon$).

\begin{theorem}
\label{thm:ABTnegative1}
For every $\Delta\ge d\ge 2$ and $g\ge 3$, there exists a positive constant
$c=c(d,\Delta,g)>0$ such that every graph $G$ of order $n$, of minimum degree $d$,
maximum degree at most $\Delta$, and with odd girth at least $g$ has
at least $cn$ eigenvalues that are larger than
$2\sqrt{d-1}\bigl(1-\varepsilon\bigr)$
and has at least\/ $cn$ eigenvalues that are smaller than
$-2\sqrt{d-1}\bigl(1-\varepsilon\bigr)$,
where $\varepsilon = \bigl(\frac{2\pi}{g}\bigr)^2 + O(g^{-3})$.
\end{theorem}

\begin{proof}
(Sketch) The proof is essentially the same as the proof of Theorem \ref{thm:main}(b),
where we take $r=\lfloor\tfrac{1}{2}g\rfloor-1$ and apply the estimate of
Theorem \ref{thm:convergence}. The assumption that the odd girth is more than $2r+1$
shows that the $r$-balls in $G$ contain no cycles of odd length. In particular, they are bipartite and hence their spectrum is symmetric with respect to 0.
Thus, knowing that the spectral radius $\l$ is large, we conclude
that the smallest eigenvalue $-\l$ is large in absolute value.
Now, we can use the interlacing theorem for the smallest eigenvalues of $G$
compared to the eigenvalues of the induced subgraph of $G$ consisting of
disjoint $r$-balls around $\lceil cn \rceil$ vertices that are $(2r+2)$-apart.
\end{proof}

The generalized version of Theorem \ref{thm:ABTnegative1} holds as well.
The proof is the same, except that we do not provide
an explicit estimate on $\varepsilon$ in terms of the odd girth.

\begin{theorem}
\label{thm:ABTnegative2}
Let\/ $D$ be a multipartite degree matrix, and let\/ $\rho_D = \rho(T_D)$.
For every $\varepsilon>0$ and every positive integer $\Delta$, there exists an integer
$g=g(D,\varepsilon)$ and a positive constant
$c=c(D,\Delta,\varepsilon)>0$ such that every graph $G$ of order $n$, of maximum
degree at most $\Delta$, with subdegree matrix $D$, and with odd girth at least $g$ has
at least $cn$ eigenvalues that are smaller than $-\rho_D+\varepsilon$.
\end{theorem}


\end{document}